\documentclass{amsproc}

\usepackage{tikz}
\usepackage[all]{xy}

\newtheorem{theorem}{Theorem}[section]

\theoremstyle{definition}
\newtheorem{definition}[theorem]{Definition}

\newtheorem{prop}[theorem]{Proposition}
\newtheorem{conj}[theorem]{Conjecture}

\theoremstyle{remark}

\numberwithin{equation}{section}

\newcommand{\ic}{\mathbb C}
\newcommand{\iz}{\mathbb Z}
\newcommand{\iq}{\mathbb Q}
\newcommand{\ndiv}{\hspace{-4pt}\not|\hspace{2pt}}

\DeclareMathOperator{\hhh}{H}

\DeclareMathOperator{\ttt}{T}
\DeclareMathOperator{\mmm}{M}
\DeclareMathOperator{\euler}{e}

\DeclareMathOperator{\hm}{Hom}
\DeclareMathOperator{\nd}{End}
\DeclareMathOperator{\gl}{GL}
\DeclareMathOperator{\hilb}{Hilb}
\DeclareMathOperator{\dt}{DT}

\begin{document}

\title[Local and relative BPS state counts]{Local and relative BPS state counts for del Pezzo surfaces}


\author[M. van Garrel]{Michel van Garrel}
\address{KIAS, 85 Hoegiro Dongdaemun-gu, Seoul 130-722, Republic of Korea}
\email{vangarrel@kias.re.kr}




\begin{abstract} Relative BPS state counts for log Calabi-Yau surface pairs were introduced by Gross-Pandharipande-Siebert. We describe how in the case of del Pezzo surfaces they are linearly related to local BPS state counts by means of generalized Donaldson-Thomas invariants of loop quivers.
\end{abstract}

\maketitle

\section{Introduction}

We describe a fascinating relationship between local and relative BPS state counts of del Pezzo surfaces, as established by the author in \cite{GWZ}. These $A$-model local and relative BPS numbers are related via a linear transformation. This correspondence is intriguing because the coefficients of the corresponding matrix are generalized Donaldson-Thomas invariants of loop quivers, which are $B$-model invariants. We as of now do not know how to explain this phenomenon. In \cite{GWZ}, these entries were computed and found to coincide with the calculation of said Donaldson-Thomas invariants as carried out by Reineke in \cite{Re,Rec}. Our description is purely mathematical, we however expect that it fits into a natural physics setting. Let us note that to our knowledge, there is as of now no mention of relative BPS state counts in the physics literature. 

\section{Local BPS state counts}
\label{sec:local}

Let $S$ be a smooth del Pezzo surface. Local BPS state counts, as defined in definition \ref{localbps} below, are the $A$-model invariants of local mirror symmetry, which was developed in \cite{CKYZ}. See also \cite{KM} of a description in terms of Yukawa couplings. Denote by $D$ a smooth effective anticanonical divisor on $S$, by $K_S$ the total space of the canonical bundle $\mathcal{O}_S(-D)$, and let $\beta\in\hhh_{2}(S,\iz)$. Denote moreover by $\overline{M}_{0,0}(S,\beta)$, resp. by $\overline{M}_{0,1}(S,\beta)$, the moduli stack of genus 0 stable maps
$$
f:C\to S
$$
with no, resp. one, marked point and such that $f_{*}([C])=\beta$. We have the forgetful morphism
$$
\pi:\overline{M}_{0,1}(S,\beta)\to\overline{M}_{0,0}(S,\beta)
$$
and the evaluation map
$$
ev:\overline{M}_{0,1}(S,\beta)\to S.
$$
Consider the obstruction bundle $R^{1}\pi_{*}ev^{*}K_{S}$. Its fibre over a stable map $f:C\to S$ is $\hhh^{1}(C,f^{*}K_{S})$. Denote by $\euler$ the Euler class and by $[\, ]^{vir}$ the virtual fundamental class.

\begin{definition}
\label{localgwdef}
The \emph{genus 0 degree $\beta$ local Gromov-Witten invariant of $S$} is
$$
I_{K_{S}}(\beta):=\int_{[\overline{M}_{0,0}(S,\beta)]^{vir}}\euler\left(R^{1}\pi_{*}ev^{*}K_{S}\right)\in\mathbb{Q}.
$$
\end{definition}

We proceed with the definition of the local BPS numbers.

\begin{definition}[See Gopakumar-Vafa in \cite{GVI, GVII} and Bryan-Pandharipande in \cite{BPB}]
\label{localbps}
Assume that $\beta$ is primitive and let $d\geq 1$. The \emph{local BPS state counts}  $n_{d\beta}\in\iq$ are defined through the following equality of generating functions.
$$
\sum_{l=1}^{\infty}I_{K_{S}}(l\beta)\, q^{l}=\sum_{d=1}^{\infty}n_{d\beta}\sum_{k=1}^{\infty}\frac{1}{k^{3}}\, q^{dk}.
$$
\end{definition}

That the BPS state counts (of any genus and for any Calabi-Yau threefold) are integers is a conjecture stated by Bryan-Pandharipande in \cite{BPB} and attributed to Gopakumar-Vafa. In genus 0, the conjecture is proven in the following case relevant to us.

\begin{theorem}[Peng in \cite{Pe}]
\label{localint}
If $S$ is a toric del Pezzo surface, then the BPS numbers $n_{d\beta}$ of $K_S$ are integers $\forall d\geq 1$.
\end{theorem}

More generally for toric Calabi-Yau threefolds, the analogous result was proven by Konishi in \cite{Ko}.

\section{Relative BPS state counts}
\label{sec:rel}

Relative BPS numbers are defined for the geometry of log Calabi-Yau surface pairs:

\begin{definition}[See section 6 of \cite{GPS}]
\label{logCY}
Let $X$ be a smooth surface, let $D\subseteq X$ be a smooth divisor and let $\gamma\in\hhh_2(S,\iz)$ be non-zero. If
$$
D\cdot\gamma = c_1(S)\cdot\gamma,
$$
then $(X,D)$ is said to be \emph{log Calabi-Yau with respect to $\gamma$}. Moreover, $(S,D)$ is said to be simply \emph{log Calabi-Yau} provided that the above equation holds for all $\gamma$, 
\end{definition}

The pairs $(S,D)$ consisting of a del Pezzo surface and a smooth anticanonical divisor form a class of examples of log Calabi-Yau surface pairs. We will proceed to defining the relative BPS state counts for this class of examples, though we note that the definition in \cite{GPS} is formulated for all log Calabi-Yau surface pairs. Their definition is entirely parallel to the definition of the local BPS state counts of section \ref{sec:local}.

We are concerned with relative stable maps, that is in addition to considering stable maps, we also prescribe the tangencies of how the maps meet $D$. Let $\beta\in\hhh_{2}(S,\iz)$ be the class of a curve and set $w=D\cdot\beta$ to be the total intersection multiplicity of $\beta$ with $D$. Denote by $\overline{M}(S/D,w)$ the moduli space of genus 0 degree $\beta$ relative stable maps meeting $D$ in one point of maximal tangency $w$. Then $\overline{M}(S/D,w)$ is of virtual dimension 0. Indeed, a generic stable map of degree $\beta$ meets $D$ in $w$ points. The moduli space of stable maps of degree $\beta$ is of virtual dimension
$
\int_\beta c_1(\ttt_S) + (\dim S -3)(1 - g) =  w-1.
$ 
Identifying two points of intersection cuts it by one. Repeating this process for all but one of the $w$ intersection points drops the virtual dimension to 0.

\begin{definition}
The \emph{genus 0 degree $\beta$ relative Gromov-Witten invariant of maximal tangency} is
$$
N_{S}[w]:=\int_{[\overline{M}(S/E,w)]^{vir}}1\in\mathbb{Q}.
$$
\end{definition}
In the notation, $\beta$ is hidden for simplicity. Relative BPS numbers are defined by extracting multiple cover contributions, under the idealized assumption that all embedded curves are rigid (and that there are only finitely many in each degree). Let $\iota:P\to S$  be such a rigid element\footnote{see \cite{GPS} for precise definitions.} of $\overline{M}(S/E,w)$ and denote by $M_{P}[k]$, for $k\geq1$, the contribution of $k$-fold multiple covers of $P$ to $N_{S}[kw]$. These numbers are calculated to be as follows:
\begin{prop}[Proposition $6.1$ in \cite{GPS}]
\[
M_{P}[k]=\frac{1}{k^{2}}\binom{k(w-1)-1}{k-1}.
\]
\end{prop}

This leads to the following definition.

\begin{definition}[Paragraph 6.3 in \cite{GPS}]
\label{def-rel-bps}
For $d\geq1$, the \emph{relative BPS state counts} $n_{S}[dw]\in\iq$ of class $d\beta$ are
defined as the unique numbers making the following equation true:
\begin{equation}
\sum_{l=1}^{\infty}N_{S}[lw]\, q^{l}=\sum_{d=1}^{\infty}n_{S}[dw]\sum_{k=1}^{\infty}\frac{1}{k^{2}}\binom{k(dw-1)-1}{k-1}\, q^{dk}.\label{eq:defn_rel_bps_gen_fn}
\end{equation}
\end{definition}

\begin{conj}[Conjecture $6.2$ in \cite{GPS}]
\label{conjecture_GPS}
For $\beta\in\hhh_{2}(S,\iz)$ an effective curve class,
$w=\beta\cdot E$ and $d\geq1$ as above,
\[
n_{S}[dw]\in\iz.
\]
\end{conj}

\begin{theorem}[Corollary 10 in \cite{GWZ}]
\label{cor}
Assume the same notation as in conjecture \ref{conjecture_GPS} and suppose furthermore that $S$ is toric. Then
\[
n_{S}[dw]\in\iz.
\]
for all $d\geq 1$.
\end{theorem}

Note that conjecture 6.2 of \cite{GPS} is stated for any Calabi-Yau surface pair.

\section{Generalized Donaldson-Thomas invariants of loop quivers}

We discuss the definition and state the explicit computation of the generalized Donaldson-Thomas invariants of loop quivers. Their definition is motivated by the framework of Kontsevich-Soibelman of \cite{KS}, and they were studied by Reineke in \cite{Re}. Their calculation as stated in \cite{Re} and reproduced in theorem \ref{thm:Re} below is a special case of a result by Reineke from \cite{Rec}.

Let $m\geq 1$, which will be fixed. The $m$-loop quiver consists of one vertex and $m$ loops. The framed $m$-loop quiver $L_m$ has an additional vertex with an arrow connecting the two vertices:

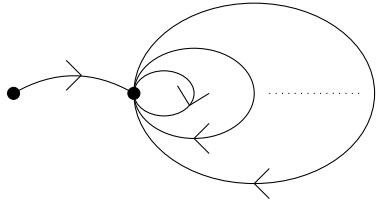
\begin{center}
\begin{tikzpicture}[scale=2]
\draw [fill] (-4.2,0) circle [radius=0.04];
\draw [fill] (-3.4,0) circle [radius=0.04];
\draw (-4.2,0) to [out=30,in=150] (-3.4,0);
\draw (-3.75,0.12) -- (-3.85,0.22);
\draw (-3.75,0.12) -- (-3.85,0.02);
\draw (-2.6,0) circle [x radius=0.8, y radius=0.6];
\draw (-2.6,-0.6) -- (-2.5,-0.5);
\draw (-2.6,-0.6) -- (-2.5,-0.7);
\draw (-3,0) circle [x radius=0.4, y radius=0.3];
\draw (-3,-0.3) -- (-2.9,-0.2);
\draw (-3,-0.3) -- (-2.9,-0.4);
\draw (-3.2,0) circle [x radius=0.2, y radius=0.15];
\draw (-3.03,-0.08) -- (-2.9,0);
\draw (-3.03,-0.08) -- (-3.11,0.05);
\draw [dotted] (-2.5,0) -- (-1.9,0);
\end{tikzpicture}
\end{center}
For $n\geq 0$, $\ic$-representations of $L_m$ of dimension $(1,n)$ consist of the data
$$
(V_0,V,\alpha_0,\alpha_1,\dots,\alpha_n),
$$
where $V_0$ and $V$ are vector spaces of dimension 1, resp. $n$, where $\alpha_0\in\hm(V_0,V)$ and where $\alpha_i\in\nd(V)$ for $i=1,\cdots,n$. A morphism $(\gamma_0,\gamma)$ between two representations is a commutative diagram
$$
\xymatrix{
V_0 \ar[r]^{\alpha_0} \ar[d]_{\gamma_0} & V \ar@(ur,rd)^{\alpha_i} \ar[d]^\gamma \\
V_0' \ar[r]^{\alpha'_0} & V' \ar@(ur,rd)^{\alpha'_i}
}
$$
such that $\gamma\circ\alpha_i=\alpha'_i\circ\gamma$. A morphism $(\gamma_0,\gamma)$ is an isomorphism if both $\gamma_0$ and $\gamma$ are. The space of all representations of $L_m$ up to isomorphism is parametrized by
$$
\ic^n\oplus\mmm_n(\ic)^{\oplus m},
$$
where $v\in\ic^n$ corresponds to $\alpha_0$ and a $m$-tuple of $n\times n$ matrices to the $\alpha_i$.

Denote by $\ic\langle x_1,\dots,x_m\rangle$ the free $\ic$-algebra on $m$ elements and let $(\phi_i)\in\mmm_n(\ic)^{\oplus m}$. Then $(\phi_i)$ determines a representation of $\ic\langle x_1,\dots,x_m\rangle$ on $\ic^m$ via $x_i\mapsto\phi_i$. Moreover, $v\in\ic^m$ is said to be \emph{cyclic} for such a representation if its image generates all of $\ic^m$, i.e. if
$$
\ic\langle\phi_1,\dots,\phi_m\rangle\, v = \ic^n.
$$
The open subset of \emph{stable} representations
$$
U \subseteq \ic^n\oplus\mmm_n(\ic)^{\oplus m}
$$
consists of those $(v,\phi_i)$ such that $v$ is cyclic for $(\phi_i)$. Consider the action of $\gl_n(\ic)$ on $U$ given by:
$$
g\cdot (v,\phi_i) = (gv,g\phi_ig^{-1}).
$$
There is a geometric quotient of this action, called the \emph{noncommutative Hilbert scheme for} $\ic\langle x_1,\dots,x_m\rangle$, and denoted by $\hilb_n^{(m)}$. We package the Euler characteristics of these spaces into a generating function:
$$
F(t):= \sum_{n\geq 0} \chi\left(\hilb_n^{(m)}\right) t^n \in \iz[[t]].
$$
Note that $F(0)=1$, so that $F(t)$ admits a product expansion.

\begin{definition}[Definition 3.1 in \cite{Re}, after \cite{KS}]
For $n\geq 1$, the \emph{generalized Donaldson-Thomas invariants  $\dt_n^{(m)}\in\iq$ of the $m$-loop quiver $L_m$} are defined by means of the following product expansion:
$$
F((-1)^{m-1}t) = \prod_{n\geq 1}(1-t^n)^{-(-1)^{(m-1)n}n\dt_n^{(m)}}.
$$
\end{definition}

Recall that the M\"obius function $\mu$, for $n\geq 1$, is defined as:
$$
\mu(n) = \begin{cases} \; \; \,1 & \text{if } n \text{ is square-free with an even number of prime factors, } \\ -1 & \text{if } n \text{ is square-free with an odd number of prime factors, } \\ \; \; \, 0 & \text{if } n \text{ is not square-free}.
\end{cases}
$$

\begin{theorem}[Reineke, theorem 3.2 in \cite{Re}]
\label{thm:Re}
$\dt_n^{(m)}\in\mathbb{N}$ and
$$
\dt_n^{(m)} = \frac{1}{n^2} \sum_{d|n} \mu\left(\frac{n}{d}\right) (-1)^{(m-1)(n-d)}\binom{mn-1}{d-1}.
$$
\end{theorem}

\section{The correspondence}

Recall that $S$ stands for a del Pezzo surface, $D$  for an anticanonical divisor on it, $\beta$ for a primitive curve class and that we set $w= D \cdot\beta$. We define an infinite-dimensional matrix $C$ with entries generalized Donaldson-Thomas invariants of loop quivers. If $t|s$, set
$$
C_{st}:=\dt_{s/t}^{(tw-1)}.
$$
If $t\ndiv s$, set $C_{st}=0$. Since lower triangular, each row of $C$ has only a finite number of non-zero entries. Hence applying to $C$ an infinite-dimensional vector does not yield convergence issues. The diagonal entries of $C$ are 1, so that $\det(C)=1$ and its inverse is integer-valued (by Cramer's rule).

We come to the interplay between the local and relative BPS state counts of $S$ that were introduced in sections \ref{sec:local} and \ref{sec:rel}. The following theorem states that the matrix $C$ provides a linear invertible relationship of the relative and local BPS numbers of $S$. Note that theorem \ref{cor}, the integrality of relative BPS state counts of toric del Pezzo surfaces, follows from the integrality of $C^{-1}$ and the theorems \ref{localint} and \ref{thm}.

\begin{theorem}[Lemma 12 in \cite{GWZ}]
\label{thm}
$$
C\cdot\left[n_{S}[dw]\right]_{d\geq 1}=\left[(-1)^{dw+1}\, dw\, n_{d\beta}\right]_{d\geq 1}.
$$
\end{theorem}

Why relative and local BPS state counts of del Pezzo surfaces should be related via Donaldson-Thomas invariants of loop quivers remains unclear.


\bibliographystyle{amsplain}

\end{document}